\documentclass[12pt,a4paper]{amsart}
\usepackage{amssymb,amsmath,amsthm}
\usepackage{hyperref}
 \usepackage[color=green!40]{todonotes}

\usepackage{lineno}

\def \bsk {\bigskip}

\def \nin {\noindent}

\newcommand\cF{{\mathcal F}}

\newcommand\cT{{\mathcal T}}

\newcommand{\ex}{\mathop{}\!\mathrm{ex}}

\newcommand{\exb}{\mathop{}\!\mathrm{ex}_\mathrm{b}}

\newcommand{\exc}{\mathop{}\!\mathrm{ex}_\mathrm{c}}

\newcommand{\exbc}{\mathop{}\!\mathrm{ex}_{\mathrm{b,c}}}

\makeatletter
\newtheorem*{rep@theorem}{\rep@title}
\newcommand{\newreptheorem}[2]{%
\newenvironment{rep#1}[1]{%
 \def\rep@title{#2 \ref{##1}}%
 \begin{rep@theorem}}%
 {\end{rep@theorem}}}
\makeatother

\theoremstyle{plain}
\newtheorem{theorem}{Theorem}
\newreptheorem{theorem}{Theorem}

\newtheorem{proposition}[theorem]{Proposition}

\newtheorem{problem}[theorem]{Problem}
\theoremstyle{definition}

\newcommand\cref[1]{Corollary~\ref{cor:#1}}

\title{Bipartite Tur\'an number of trees}

\author{Yair Caro}
\address{Department of Mathematics, University of Haifa-Oranim, Israel}
\email{yacaro@kvgeva.org.il}

\author{Bal\'azs Patk\'os}
\address{HUN-REN Alfr\'ed R\'enyi Institute of Mathematics}
\email{patkos@renyi.hu}

\author{Zsolt Tuza}
\address{HUN-REN Alfr\'ed R\'enyi Institute of Mathematics and University of Pannonia}
\email{tuza.zsolt@mik.uni-pannon.hu}

\date{}

\begin{document}


\begin{abstract}
We start a systematic investigation concerning bipartite Tur\'an number for trees.
For a graph $F$ and integers $1 \leq a \leq b$ we define:

$(i)$\quad 
$\exb(a, b, F)$ is the largest number of edges that an $F$-free bipartite graph can have with part sizes $a$ and $b$. We write $\exb(n, F)$ for $\exb(n, n, F)$.

$(ii)$\quad 
$\exbc(a, b, F)$ is the largest number of edges that an $F$-free connected, bipartite graph can have with part sizes $a$ and $b$. We write $\exbc(n, F)$ for $\exbc(n, n, F)$. 

Both definitions are similar for a family $\cF$ of graphs.

We prove general lower bounds depending on the maximum degree of $F$, 
as well as on the cardinalities of the two vertex classes of $F$.

We derive upper and lower bounds for $\exb(n,F)$ in terms of $\ex(2n,F)$ and $\ex(n, F)$, the corresponding classical (not bipartite) Tur\'an numbers.  

We solve both problems for various classes of graphs, including all trees up to six vertices for any $n$, for double stars $D_{s ,t}$ if $a \geq f(s,t )$, for some families of spiders, and more.

We use these results to supply an answer to a problem raised by L. T. Yuan and X. D. Zhang [{\it Graphs and Combinatorics}, 2017] concerning $\exb( n, \cT_{k,\ell} )$, where $\cT_{k,\ell}$ is the family of all trees with vertex classes of respective cardinalities $k$ and $\ell$.

The asymptotic worst-case ratios between Tur\'an-type functions are also inverstigated.

\bigskip


\bigskip

\end{abstract}

\maketitle

\section{Introduction and results}

In this paper we initiate a systematic study of Tur\'an numbers of
 trees in the bipartite setting.
Namely, we raise the analogous problem with restriction to the bipartite universe, i.e., where the extremal graphs are subgraphs
 of the complete bipartite graphs $K_{a ,b}$, mostly of $K_{n,n}$.
Moreover, we also consider the situation where the extremal graph
 is required to be connected.
Part of our motivation is to complement the recent research
 concerning connected Tur\'an numbers of trees without the bipartite
 restriction started in \cite{CPT} and continued in \cite{JLS}.
Another part of the motivation comes from successfully extending and answering extremal results explicitly raised in
\cite[Question 1.5]{YZ} concerning $\cT_{k,l}$, the family of all trees $T(k, l)$ with vertex classes of respective cardinalities $k$ and $l$.

Sporadic results concerning both problems, in the bipartite setting or the non-bipartite connected version already appeared; we refer to \cite{CWYZ}, \cite{GyRS}, \cite{K}, \cite{LN}, \cite{YZ}, \cite{YZ2}.
Later we shall say more about these papers, as well as on the progress we made concerning some questions raised in \cite{YZ}.  
 
Before proceeding with the main results of this paper, we introduce a necessary definition and notation which will be used in the sequel.

\bsk

\noindent
\textbf{Definition.}
Let $F$ be a given ``forbidden'' graph.
\begin{itemize}
    \item $\ex(n,F)$ is the largest number of edges that an $n$-vertex graph can have without any subgraph isomorphic to $F$. (Such graphs are called \textit{$F$-free}).
    \item $\ex_c(n,F)$ is the largest number of edges that an $n$-vertex \textit{connected} $F$-free graph can have.
    \item $\exb(a,b,F)$ is the largest number of edges that an $F$-free \textit{bipartite} graph can have with part sizes $a$ and $b$. We write $\exb(n,F)$ for $\exb(n,n,F)$. 
    \item $\exbc(a,b,F)$ is the largest number of edges that an $F$-free \textit{connected, bipartite} graph can have with part sizes $a$ and $b$. We write $\exbc(n,F)$ for $\exbc(n,n,F)$. 
\end{itemize} 
All parameters can be defined analogously when we forbid copies of all graphs of a family $\cF$. Then we use the notation $\ex(n,\cF)$, $\exc(n,\cF)$, $\exb(n,\cF)$, $\exbc(n,\cF)$.

\bsk

\noindent
\textbf{Notation}.
We use the standard notation $\delta=\delta(G)$, $\Delta=\Delta(G)$, and $e(G)$, respectively, for the minimum degree, maximum degree, and the number of edges of a given graph $G$.
The vertex-disjoint union of two
graphs $G$ and $H$ is denoted by $G\cup H$.
For a graph $H$ and a positive integer $k$, $kH$ denotes the pairwise vertex-disjoint union of $k$ copies of $H$.
We denote by
$S_k$ the star with $k$ leaves; $P_k$ and $C_k$ denote the path and the cycle on $k$ vertices, respectively. The complete bipartite graph with parts of size $a$ and $b$ is denoted by $K_{a,b}$.
We write
$D_{a,b}$ to denote the \textit{double star} on $a+b+2$ vertices such that the two non-leaf vertices have degree $a+1$ and $b+1$. The \textit{spider} $S_{a_1,a_2,\dots, a_j}$ with $j\ge 3$ is the graph obtained from $j$ paths with $a_1,a_2,\dots,a_j$ edges by identifying one endpoint of all paths. So $S_{a_1,a_2,\dots,a_j}$ has $1+\sum_{i=1}^ja_i$ vertices and maximum degree $j$. The only vertex of degree at least 3 is the \textit{center} of the spider, and the maximal paths starting at the center are its \textit{legs}. If a spider has $h$ legs of length $\ell$, and other legs of length $a_1,\dots,a_j$, then we use the notation $S_{a_1,\dots,a_j,h*\ell}$. Finally, $P_{r,s,t}$ denotes the tree on $r+s+t+3$ vertices $u,u_1,\dots,u_r$, $v,v_1,\dots,v_s$, $w,w_1\dots,w_t$ with $uvw$ being a path and $zz_i$ forming an edge for all $z\in\{u,v,w\}$ and all feasible values of $i$; hence, $P_{r,s,t}$ is a caterpillar whose spine is the path $P_3$.

\bsk

The following list covers the references relevant to this paper, including the sporadic publications of recent years concerning the connected Tur\'an numbers in the general case, and Tur\'an numbers in the bipartite setting.  
For a general survey on Tur\'an numbers see \cite{FS}, though in some sections it is a bit outdated.
The exact determination of Tur\'an numbers for paths, also  including the connected version, was done long ago in \cite{GyRS} and \cite{K}.
Results concerning cycles in bipartite graphs were considered several times; see \cite{A}, \cite{Gy}, \cite{J}, \cite{JM}, \cite{LN}.
Tur\'an numbers of linear forests in bipartite graphs were considered in \cite{CLZ}, \cite{CWYZ}, \cite{YZ2}.
Connected Tur\'an number for trees is systematically studied in \cite{CPT}, \cite{JLS}].
The famous Erd\H os--S\'os conjecture on trees appeared first in \cite{E}, and a weak version of it for the bipartite case appeared in \cite{YZ}.

\bsk

\subsection{Our contribution and main results of this paper}

\ 

\bsk

\noindent
\textbf{Preliminary assumptions.}
Throughout this paper, without further mention, we assume that the forbidden graph $F$ under consideration is \textit{bipartite} and has \textit{no isolated vertices}.
In order to avoid trivialities, in the study of $\exb(a,b,F)$ and $\exbc(a,b,F)$ with $a\leq b$ (possibly $a=b=n$), by default we also assume:
 \begin{itemize}
     \item $|V(F)|\leq a+b$;
     \item if $F$ has part sizes $p$ and $q$, then $p\leq a$ and $q\leq b$.
 \end{itemize}
Moreover, in the context of $\exbc(n,F)$ it will be assumed that $F$ is none of the graphs $K_2$, $2K_2$, $P_3$, $P_3\cup K_2$, $P_4$.

\bsk

The explanation for the exclusions with respect to $\exbc$ is the following necessary and sufficient condition for $\exbc(n,F)$ to exist for all $n$.

\begin{proposition} \label{p:exist}
    For every $n$, there exists a spanning connected subgraph of $K_{n,n}$ not containing a graph $F$ if and only if
    $F\notin \{ K_2, 2K_2, P_3, P_3\cup K_2, P_4 \}$.
\end{proposition}

Some basic general inequalities relating bipartite Tur\'an numbers to the classical Tur\'an numbers are given in the next assertion.

\begin{theorem}   \label{p:old-1-3}
For any graph $F$ and any $n$,
$  \textstyle
\ex(n, F) < \exb(n, F) \le \ex(2n,F)
$,
and $\exb(n, F) \geq (\frac{1}{2}+\frac{1}{4n-2})\ex(2n,F)$.
\end{theorem}

 The next proposition collects several lower bounds with regards to the maximum degree $\Delta(F)$.

    \begin{proposition}   \label{p:lower}
    In terms of the maximum degree $\Delta(F)$, the following lower bounds are valid.
\begin{enumerate}
    \item If $\Delta(F)\ge 2$, then $\exb(n,F)\ge (\Delta(F)-1)n$.
    \item If $\Delta(F)=1$ and $F\neq K_2$, then $\exb(n,F)\ge n$.
    \item If $\Delta(F)\ge 3$, then $\exbc(n,F)\ge (\Delta(F)-1)n$.
    \item If $F\notin \{ K_2, 2K_2, P_3, P_3\cup K_2, P_4 \}$, then $\exbc(n,F)\ge 2n-1$.
\end{enumerate}
        In particular, $\exb(n,F)\ge n$ for every $F$ with more than one edge.
\end{proposition}

The next proposition collects further types of lower bounds for both $\exb(n, F)$ and $\exbc(n, F)$, which will be used in the sequel for graphs $F$ that are connected.

\begin{proposition}\label{lower}
    In terms of part sizes, the following lower bounds are valid.
    \begin{enumerate}
        \item 
        If $s$ is the smaller part size of a connected bipartite $F$, then $\exb(a,b,F)\ge (s-1)(a-s+1+b-s+1)$. In particular, $\exb(n,F)\ge 2(s-1)(n-s+1)$.
        \item 
        If $L$ is the larger part size of a connected bipartite $F$, then $\exb(n,F)\ge \lfloor\frac{n}{L-1}\rfloor\cdot (L-1)^2$. Similarly, $\exb(n,F)\ge \lfloor\frac{2n}{|V(F)|-1}\rfloor\cdot \lfloor\frac{|V(F)|-1}{2}\rfloor \lceil\frac{|V(F)|-1}{2}\rceil$.
    \end{enumerate}
\end{proposition}

Related to maximum degree, the case of stars and the ``star+leaf'' graphs can be settled easily for $a=b=n$.

\begin{proposition}   \label{p:star}
    For every $d\geq 2$ and every $n\geq |V(F)|$ with
    $F\in\{K_{1,d+1},S_{2,d*1}\}$,
$$
    \exb(n,K_{1,d+1}) = \exbc(n,K_{1,d+1}) = 
    \exb(n,S_{2,d*1}) = 
    \exbc(n,S_{2,d*1}) = dn\,.
$$
    For the case of $d=1$, i.e.\ $K_{1,2}\cong P_3$ and $S_{2,1}\cong P_4$, we have $\exb(n,P_3) = n$ and $\exb(n,P_4) = 2n-2$, while $\exbc(n,P_3)$ and $\exbc(n,P_4)$ are not defined.
\end{proposition}

\bsk

\subsection{Small graphs}
\label{ss:small}

\ 

\bsk

This subsection is a warm-up to make the reader familiar with the extremal functions under consideration.
In two theorems we summarize bipartite
Tur\'an numbers and their connected versions for all those trees on at most six vertices which have not been considered so far.
Stars and ``star+leaf'' graphs are missing from both lists, as these simple types of graphs are settled in Proposition \ref{p:star} above.
Also, without the connectivity assumption, the case of paths (of any length) was solved by Gyárfás, Rousseau, and Schelp \cite{GyRS}.
In particular, $\exb(n,P_5)$ is equal to 
$2n-1$ if $n\geq3$ is odd, and $2n$ if $n\geq 4$ is even.
Moreover, $\exb(3,P_6)=6$, having its extremal graph $K_{2,3}\cup K_1$, and $\exb(n,P_6)=4n-8$ with extremal graph $2K_{2,n-2}$ for all $n\geq 4$.
Therefore, below we consider $P_5$ and $P_6$ in the second result only.

The values for trees of orders $4,5,6$ are collected in Table \ref{tab:small}.
For certain trees there are small exceptional cases, which are not listed in the table; to see them, please see the corresponding results directly.

\renewcommand{\arraystretch}{1.3}
\begin{table}[ht]
    \centering
    \begin{tabular}{ccccccccc}
        $n$ && graph && $\exb$ && $\exbc$ && reference \\
    \hline
    \hline
        4 && $P_4$ && $2n-2$ && $-$ && \cite{GyRS} \& Proposition \ref{p:exist} \\
         && $K_{1,3}$ && $2n$ && $2n$ && Proposition \ref{p:star} \\
    \hline
        5 && $P_5$ && $2n-(n\, \mathrm{mod}\, 2)$ && $2n-1$ && \cite{GyRS} \& Theorem \ref{smalltreesconnected} \\
         && $K_{1,4}$ && $3n$ && $3n$ && Proposition \ref{p:star} \\
         && $S_{2,1 ,1}$ && $2n$ && $2n$ && Proposition \ref{p:star} \\
    \hline
        6 && $P_6$ && $4n-8$ && $2n-1$ && \cite{GyRS} \& Theorem \ref{smalltreesconnected} \\
         && $K_{1,5}$ && $4n$ && $4n$ && Proposition \ref{p:star} \\
         && $S_{3,1,1}$ && $3n-c_{n \, \mathrm{mod} \, 3}$ && $2n$ && Theorems \ref{spiders} \& \ref{smalltreesconnected} \\
         && $S_{2 ,2 ,1}$ && $4n-8$ && $2n$ && Theorems \ref{smalltrees} \& \ref{smalltreesconnected} \\
         && $S_{2 ,1 ,1 ,1}$ && $3n$ && $3n$ && Proposition \ref{p:star} \\
         && $D_{2,2}$ && $4n-8$ && $4n-9$ && Theorems \ref{smalltrees} \& \ref{smalltreesconnected} \\
    \hline
         &&  &&  &&  &&  \\
    \end{tabular}
    \caption{General values of $\exb(n,T)$ and $\exbc(n,T)$ for $4\leq |V(T)|\leq 6$, disregarding possible small exceptions.}
    \label{tab:small}
\end{table}

\begin{theorem}\label{smalltrees} \
\begin{enumerate}
    \item 
    $\exb(n,S_{3,1,1})\leq 3n$, with equality if and only if $n$ is a multiple of $3$.
    In that case the unique extremal graph is $\frac{n}{3}K_{3,3}$.
    \item 
    $\exb(n,D_{2,2}) = 4n -8$ for $n \ge 5$, with the unique extremal graph $2K_{2,n-2}$ if $n\geq 6$.
    Moreover, $\exb(3,D_{2,2}) = 6 = 4n -6$ for $n = 3$ and $\exb(4,D_{2,2}) = 9 = 4n -7$ for $n = 4$.
    \item 
    $\exb(n,S_{2,2,1})=4n-8$ if $n\geq 4$ and $\exb(3,S_{2,2,1})=6=4n-6$.
\end{enumerate}
\end{theorem}

Some bipartite connected Tur\'an numbers are determined already in the above Theorem \ref{smalltrees} and Proposition \ref{p:star}. 
For the other trees up to six vertices, we state the next result.
Note that $F=P_3$ and $F=P_4$ are infeasible with respect to $\exbc$, due to Proposition \ref{p:exist}.

\begin{theorem}\label{smalltreesconnected}
\

\begin{enumerate}
    \item 
    $\exbc(n,P_5) =
    \exbc(n,P_6) = 2n-1$.
    \item 
    $\exbc(n,D_{2,2})=4n-9$ for all $n\geq 6$ and $\exbc(n,D_{2,2})=\exb(n,D_{2,2})$ for $3\leq n\leq 5$.
    \item 
    $\exbc(n,S_{3,1,1})=2n$.
    \item 
    $\exbc(n,S_{2,2,1})=2n$.
\end{enumerate}
    
\end{theorem}

\bsk

\subsection{Some spiders and other classes of trees}

\ 

\bsk

Here we state tight results and general upper bounds on some infinite classes of trees.
In these theorems the general host graph $K_{a,b}$ is considered, i.e., not restricted to the balanced case of $K_{n,n}$.
The first result also puts Proposition \ref{p:star} into a more general setting.

\begin{theorem}\label{spiders} \
    \begin{enumerate}
        \item 
        For $d\ge 2$ and $a\ge d+1$, $\exb(a,b,S_{2,d*1})=\max\{da,d(a-1)+(b-a+1)\}$.
        \item 
        For $d\ge 1$ we have $\exb(n,S_{3,d*1})=(d+1)n - c_{d,n}$ (where $c_{d,n}>0$ occurs when $n$ is not a multiple of $d+1$).
        \item 
        For every $d\geq 2$, $\exbc(n,S_{3,d*1})=dn$.
    \end{enumerate}
\end{theorem}

Note that for $d=1$ we have $S_{3,1*1}=P_5$, a graph discussed above already.

\begin{theorem}\label{classes} \
    \begin{enumerate}
        \item 
        $\exb(a,b,D_{s,t})=s(a+b-2s)$ if $a,b\ge 4(s+1)^3$ and $2s>t\ge s$. Moreover, the unique extremal graph is $K_{s,a-s}\cup K_{s,b-s}$. 
        \item
        $\exbc(a,b,D_{s,t})=s(a+b-2s)-1$ if $a,b\ge 4(s+1)^3$ and $2s>t\ge s$.
        
        \item 
        If $2s\le t$, then $\exb(a,b,D_{s,t})\le t(a+b)/2$ and the inequality is sharp if $a=b=n$ is divisible by $t$.
        \item 
        For any $\ell\ge 2$, $\exb(a,b,S_{2\ell,(\ell+1)*1})\le \ell (a+b)$ and the inequality is sharp if $a=b=n$ is divisible by $2\ell$.
        \item 
        If $s< \max\{r,t\}=:m$, then $\exb(a,b,P_{r,s,t})\le m\cdot (a+b)$.
    \end{enumerate}
\end{theorem}

Let $\cT_{k,\ell}$ denote the set of trees that have bipartition part sizes $k$ and $\ell$. Yuan and Zhang \cite{YZ} started to investigate $\exb(a,b,\cT_{k,\ell})$ and determined $\exb(a,b,\cT_{k,2})$ and $\exb(a,b,\cT_{3,3})$ for all values of $a$ and $b$. 
The next theorem gives exact values of $\exb(a,b,\cT_{k,\ell})$ for wide ranges of $( a ,b ,k , \ell)$,
supplying many further extremal results to a problem explicitly raised in \cite[Question 1.5]{YZ}.

\begin{theorem}   \label{t:trees-k-l}
    Suppose $k\le \ell$ and  $a\le b$.
    \begin{enumerate}
        \item 
        If $k\le \ell\le 2k-1$ and $a\ge 4k^3$, then $\exb(a,b,\cT_{k,\ell})=(k-1)(a+b-2(k-1))$ and the unique extremal graph is
        $K_{k-1,a-k+1} \cup K_{k-1,b-k+1}$.
        \item 
        If $2k\le \ell$, then $\exb(a,b,\cT_{k,\ell})\le (\ell-1)(a+b)$ and the inequality is sharp if $a = b = n$ and $n$ is
        divisible by $\ell-1$.
    \end{enumerate}
\end{theorem}

\bsk

\subsection{Worst-case ratios}

\ 

\bsk

In \cite{CPT}, we introduced the parameter $\gamma$ to denote the smallest ratio $\exc(n,T)$ and $\ex(n,T)$ can have. More precisely, let $\cT_k$ denote the set of all trees on at least $k$ vertices. For a tree $T$, we define $\gamma_T=\limsup_n\frac{2}{|T|-2}\cdot \frac{\exc(n,T)}{n}$ and then $\gamma=\lim_{k \rightarrow \infty}\inf\{\gamma_T: T\in \cT_k\}$. If the Erd\H os--S\'os conjecture holds true, then $\ex(n,T)\le \frac{|T|-2}{2}n$ for all trees, so the ratio is at most 1. In \cite{CPT}, $1/3\le \gamma \le 2/3$ was proved, and then Jiang, Liu, and Salia \cite{JLS} settled $\gamma=1/2$. In the same flavor we define 

$$\gamma_{b,T}=\limsup_n\frac{\exb(n,T)}{(|T|-2)n}, \hskip 0.5truecm \gamma_{b,c,T}=\limsup_n\frac{\exbc(n,T)}{(|T|-2)n}, \hskip 0.5truecm \gamma_{b \leftrightarrow c,T}=\limsup_n\frac{\exbc(n,T)}{\exb(n,T)}$$ 

\nin
and 

$$\gamma_{b}=\lim_{k \rightarrow \infty}\inf\{\gamma_{b,T}: T\in \cT_k\}, \gamma_{b,c}=\lim_{k \rightarrow \infty}\inf\{\gamma_{b,c,T}: T\in \cT_k\}, \gamma_{b\leftrightarrow c}=\lim_{k \rightarrow \infty}\inf\{\gamma_{b \leftrightarrow c,T}: T\in \cT_k\}.$$ 

\begin{theorem}\label{gammab}
    $\gamma_b=2/3$.
\end{theorem}

By definition $\gamma_{b,c}\le \gamma_b$ holds.
Also, clearly we have $\exb(n,T)\le \ex(2n,T)$, moreover the Erd\H os--S\'os conjecture would imply $\ex(2n,T)\le (|T|-2)n$ and so $\gamma_{b,c}\le \gamma_{b\leftrightarrow c}$ would follow. A well-known greedy embedding strategy shows $\ex(2n,T)\le 2(|T|-2)n$, hence $\gamma_{b,c}\le 2\gamma_{b\leftrightarrow c}$.

\begin{theorem}\label{gammabc}
    $\frac{2-2x_0}{3}\le \gamma_{b,c}$, where $x_0$ is the larger root of the equation $14x^2-14x+3=0$, $x_0=0.68898...$. And so $\gamma_{b,c}\ge 0.207$ and $\gamma_{b\leftrightarrow c}\geq 0.1035$. 
\end{theorem}

Proposition \ref{p:star} and the theorems of Section \ref{ss:small} will be proved in Section \ref{s:small} at the end of the paper, all the other results in the next Section \ref{s:general}.

\bsk

\section{Proofs for general tree classes}
\label{s:general}

\bsk

\begin{proof}[{\bf Proof of Proposition \ref{p:exist}}]
Consider the cycle, the double star, and the spider
$$
  C_{2n} \,, \qquad D_{n-1,n-1} \,, \qquad S_{1,(n-1)*2} \,.
$$
Each of them is a spanning connected subgraph of $K_{n,n}$.
Hence, if $\exbc(n,F)$ is undefined, then $F$ is a subgraph of all the three.
Assume that $F$ has no isolated vertices.
From $C_{2n}$ we see that every component of $F$ is a path.
Then $D_{n-1,n-1}$ implies:
\begin{itemize}
    \item there are at most two components (paths);
    \item the longest component is $K_2$ or $P_3$ or $P_4$;
    \item $P_4$ does not admit a second component.
\end{itemize}
Thus, beyond the graphs $K_2$, $2K_2$, $P_3$, $P_3\cup K_2$, $P_4$ the only possibility that remains is $P_3\cup P_3$.
But in $S_{1,(n-1)*2}$ any two copies of $P_3$ share a vertex.
This completes the proof of the ``if'' part.

For the ``only if'' part it suffices to consider $P_4$ and $P_3\cup K_2$.
Let $H$ be any connected spanning subgraph of $K_{n,n}$, $n\geq 3$.
We need to prove $P_4\subset H$ and $P_3\cup K_2\subset H$.
If $P_4\not\subset H$ held, then $H$ would be a star, but then it would not span $K_{n,n}$.
Concerning $P_3\cup K_2$, let $v$ be a highest-degree vertex of $H$, and $u_1,\dots,u_k$ its neighbors.
Here $k\geq 2$ because $H$ is connected.
If $k\geq 3$, we pick any $v'\neq v$ from the vertex class of $v$, and any edge $v'u'$ (as $H$ has no isolates).
We may assume $u'\notin\{u_1,u_2\}$, thus the three edges $vu_1,vu_2,v'u'$ form $P_3\cup K_2$.
If $k=2$, pick a third vertex $u'$ in the class of $u_1,u_2$ and its any neighbor $v'$.
Since $v\neq v'$, also here the three edges $vu_1,vu_2,v'u'$ form $P_3\cup K_2$.
\end{proof}

\begin{proof}[{\bf Proof of Theorem \ref{p:old-1-3}}]
    The upper bound $\ex(2n, F)$ is trivial, as $\exb(n,F)$ considers graphs on $2n$ vertices. 
Also, the lower bound $\ex(n, F)$ will follow once we prove the same for $\frac{1}{2}\ex(2n, F)$.
Indeed, the union of two vertex-disjoint copies of an $F$-extremal graph of order $n$ shows $\ex(2n, F) \geq 2\ex(n, F)$.

For the proof of the last inequality $\exb(n, F) \geq (\frac{1}{2}+\frac{1}{4n-2})\ex(2n,F)$,
let $G$ be any $F$-free graph of $2n$ vertices and $\ex(2n, F)$ edges.
We start with an arbitrary balanced vertex partition $A\cup B=V(G)$, and apply the local switching algorithm due to Bylka, Idzik, and Tuza \cite{BIT}.
Namely, if there exists an edge $uv\in E(G)$ with $u\in A$ and $v\in B$ such that the partition $(A',B')$ with $A'=(A\setminus\{u\})\cup\{v\}$, $B'=(B\setminus\{v\})\cup\{u\}$ generates a larger edge cut (i.e., there are more $A'$--$B'$ edges than $A$--$B$ edges) then we replace $(A,B)$ with $(A',B')$.
Repeat this step as long as such edges exist.
According to Theorem 3.2 of \cite{BIT}, the cut contains at least $\frac{n}{2n-1}\,|E(G)|$ edges when no further improving step is possible.
\end{proof}

We continue with the lower bounds which use different parameters of the forbidden graph $F$.

\begin{proof}[{\bf Proof of Proposition \ref{p:lower}}]
    The edge set of $K_{n,n}$ can be partitioned into $n$ perfect matchings, so $(\Delta(F)-1)$-regular subgraphs of $K_{n,n}$ exist.
To ensure connectivity, for $\Delta(F)-1=2$ we take a Hamiltonian cycle $C_{2n}$.
For $\Delta(F)-1>2$, we supplement $C_{2n}$ with any $\Delta(F)-3$ perfect matchings edge-disjoint from it.
An explicit way to do so is to label the vertices in the two parts of $K_{n,n}$ as $u_1,\dots,u_n$ and $v_1,\dots,v_n$ and join each $u_i$ ($i=1,\dots,n$) with $v_i,v_{i+1},\dots,v_{i+k-1}$, subscript addition taken modulo $n$.

The lower bound (4) follows from Proposition \ref{p:exist}, as the number of edges in a spanning tree of a $2n$-vertex graph is $2n-1$.
\end{proof}

\begin{proof}[{\bf Proof of Proposition \ref{lower}}]
    If $F$ is connected and its smaller part size is $s$, then clearly $K_{s-1,a-s+1}\cup K_{s-1,b-s+1}\subset K_{a,b}$ is $F$-free. This shows (1).

    If $n=x(L-1)+y$ with $0\le y<L-1$, then $xK_{L-1,L-1}\cup K_{y,y}\subseteq K_{n,n}$ is $F$-free. Also, if $x=\lfloor \frac{n}{|V(F)|-1}\rfloor$, then $2xK_{\lfloor \frac{|V(F)-1}{2}\rfloor,\lceil \frac{|V(F)-1}{2}\rceil}\subset K_{n,n}$ is $F$-free. This shows (2).
%
\end{proof}

\begin{proof}[{\bf Proof of Proposition \ref{p:star}}]
Since both $K_{1,d+1}$ and $S_{2,d*1}$ have maximum degree $d+1\geq 3$ if $d>1$, parts (1) and (3) of Proposition \ref{p:lower} imply that $dn$ is a lower bound on all the four values.
If $d=1$, then the perfect matching $nK_2$ and the graph $2K_{1,n-1}$ provide a construction for $P_3$ and $P_4$, respectively.
Tightness of the latter two is simply explained by the facts that the $P_3$-free graphs are exactly the matchings (together with any number of isolated vertices), and every connected component of a $P_4$-free subgraph of $K_{n,n}$ is a star or an isolated vertex.
Having $c$ such components, the number of edges is $2n-c$, and a single star cannot span $K_{n,n}$ if $n\geq 2$.

Concerning the upper bound $dn$ it suffices to consider $\exb(n,S_{2,d*1})$, because
it is at least as large as any of the other three values under consideration.
Assume that $G\subset K_{n,n}$ is $S_{2,d*1}$-free.
If a vertex $v$ has degree at least $d+1$, then all its neighbors have degree 1, for otherwise an $S_{2,d*1}$ with center $v$ would occur.
Thus, each component of $G$ either is a star or its maximum degree does not exceed $d$.
Consequently the average degree in $G$ is at most $d$, and the number of edges is at most $dn$.
\end{proof}

\begin{proof}[{\bf Proof of Theorem \ref{spiders}}]
To see the upper bound of (1), observe that if a vertex $v$ in a bipartite $S_{2,d*1}$-free subgraph $G$ of $K_{a,b}$ has degree at least $d+1$, then the component of $v$ in $G$ must be a star. Indeed, $v$'s neighbors form an independent set as $G$ is bipartite, and a neighbor of a neighbor would give rise to a copy of $S_{2,d*1}$. So if $m$ is the number of vertices in the smaller part of $G$ in non-star components (in particular $m\le a$), then $e(G)\le dm+(a-m+b-m-1)$, while if $a=m$, then there is no star-component of $G$ and thus $e(G)\le da$. The former expression is maximized at $m=a-1$ as $d\ge 2$. This proves the upper bound. The construction giving the $ad$ lower bound is a $d$-regular subgraph of $K_{a,a}\subseteq K_{a,b}$. The construction giving the $(a-1)d+b-a+1$ lower bound is $H\cup K_{1,b-a+1}$ where $H$ is a  $d$-regular subgraph of $K_{a-1,a-1}$.

\medskip

The construction for the lower bound of (2) is as follows: vertex-disjoint copies of $K_{d+1,d+1}$. This is indeed $S_{3,d*1}$-free, since the bipartition of $S_{3,d*1}$ has vertex-class distribution $(2,d+2)$. This is supplemented with $K_{b,b}$ if $n$ is of the form $(d+1)a+b$ with $0<b<d+1$. 

\medskip

Proof of the upper bound of (2): we prove that if the maximum degree of a bipartite graph exceeds $d+1$, then it cannot be extremal for $S_{3,d*1}$,  hence $\Delta(G) \le  d+1$ and then $e(G) \le (d+1)\cdot 2n /2 = (d+1)n$. 

Let $G$ be any $S_{3,d*1}$-free subgraph of $K_{n,n}$. Assume that a vertex $v$ has $m > d+1$ neighbors, say $N(v)=U=\{u_1,u_2,\dots,u_m\}$. We denote by $Z$ the set of vertices that have at least one neighbor in $U$.
Let $z$ be any vertex from $Z\setminus \{v\}$.
If $z$ has more than one neighbor in $U$, say both $zu_1$ and $zu_2$ are edges, then we take $v$ as the center of $S_{3,d*1}$, with long leg $u_1,z,u_2,v$ and short legs $vu_3, vu_4,\dots, vu_{d+2}$.
If $z$ has at least one neighbor say $u_0$ outside $U$, then $u_0$ can take the previous role of $u_1$ and again an $S_{3,d*1}$ centered at $v$ is found. Thus $z$ is a pendant vertex in $G$.

It follows that if a connected component $G'$ of $G$ has maximum degree bigger than $d+1$, then it is a tree (of radius 2).
In particular, if $G'$ has $p$ vertices in one vertex class and $q$ vertices in the other, then it has $p+q-1$ edges.

We now consider all the tree components in $G$ together. They contain some number of vertices in each vertex class;
let $s$ denote the larger of those two numbers. Note that $s$ is at least $d+2$, already due to $N(v)$. Then the number of edges in $G$ is at most
\[2s-1 + (d+1)(n-s) = (d+1)n - (d-1)s -1 \leq (d+1)n - (d-1)(d+2) -1.
\]
Indeed, $2s-1$ is an upper bound in the tree components, and we enumerate each of the other edges at its vertex in the vertex class where the trees together cover $s$ vertices.

So the loss compared to $(d+1)n$ is at least $d^2 + d -1 > d^2$ if $d>1$.
On the other hand our construction takes as many copies of $K_{d+1,d+1}$ as possible, and if $n=(d+1)a+b$, then we have an additional $K_{b,b}$. So the loss in this construction is $(d+1)b-b^2 = (d+1-b)b \leq (d+1)^2 /4 < d^2$ if $d>1$, therefore $G$ cannot be extremal.

\medskip

The lower bound in (3) is verified by any $d$-regular, connected, bipartite graph.
For the upper bound let $G\subset K_{n,n}$ be any graph with $dn+1$ edges.
Then both vertex classes $A,B$ contain at least one vertex of degree at least $d+1$; say, these are $u\in A$ and $v\in B$.
If $N(u)\neq B$, then consider a shortest path $P$ from $B\setminus N(u)$ to $N(u)$.
In this case an $S_{3,d*1}$ is present with center $u$, whose legs are $P$ together with $d$ leaves from $N(u)\setminus V(P)$.
The same situation occurs if $N(v)\neq A$.
And if both $N(u)=B$ and $N(v)=A$ hold, then the double star $D_{n-1,n-1}$ (whose central edge is $uv$) together with any one further edge forms an $S_{3,d*1}$; the center can be taken as either of $u$ and $v$.
\end{proof}

\begin{proof}[{\bf Proof of Theorem \ref{classes}}]
The lower bound of (1) is given by the construction $K_{s,a-s}\cup K_{s,b-s}$. For the upper bound, it is enough to consider the case $t=2s-1$. Let $G$ be an extremal bipartite graph with vertex classes $A,B$ where $|A|=a,|B|=b$.

We use the following notation:
the set of vertices of degree greater than $2s-1$ in $A$ and $B$ are denoted $A^H$, $B^H$ respectively.
The set of vertices of degree at most $s$ in $A$ and $B$ are denoted by $A^L$, $B^L$ respectively.
The set of vertices of degree in $[s+1, 2s-1]$ are denoted by $A^M, B^M$.
(Superscripts represent High / Middle / Low degrees.)
 
All edges incident with $A^H \cup B^H$ have their other end in $A^L \cup B^L$ as otherwise $G$ contains a $D_{s,2s-1}$.
For each of those edges we assign weight 1 to their endvertex in $A^L \cup B^L$
and weight 0 to their end in $A^H \cup B^H$.
For every other edge we assign weight $1/2$ at both of its ends.
An immediate consequence is that every vertex having at least one neighbor
in $A^M \cup B^M$ has weight at most $s - 1/2$, hence at least $1/2$ less than $s$.
 
The sum of weights equals the number of edges, moreover all vertices
in $A^H \cup B^H$ have weight 0. Thus, as a basic upper bound, we have
\begin{equation}\label{eq1}
    |E| \leq s(a+b - A^H - B^H) \tag{*}
\end{equation}   
with equality if and only if every vertex in $A^L \cup B^L$ has degree exactly $s$,
moreover $A^M \cup B^M$ is empty and there are no edges inside $A^L \cup B^L$.
For a more detailed computation we distinguish between two cases.

\medskip

\textsc{Case I} Both $A^H$ and $B^H$ contain at least $s$ vertices.

Then (\ref{eq1}) immediately yields the required upper bound, and also uniqueness of the extremal graph $K_{a,n-s}\cup K_{s,b-s}$ follows in this case, using the comments
given under (\ref{eq1}).

 \smallskip 
 
\textsc{Case II}  $A^H$ (or similarly $B^H$, or both) contains fewer than $s$ vertices.

Assume that $|A^H| = s  - x$ for some $x > 0$.
Then the vertices in $B^L$ have fewer than $s$ neighbors in $A^H$, therefore their
weights are at most $s-x/2$. Hence all vertices in $A^M \cup B^M \cup B^L$ have weight
at most $s-1/2$. So 
\[|E|\leq s|A^L|+(s-1/2)(|A^M| + |B^M| + |B^L|)\leq s(a+b) - 1/2 (|A^M| + |B^M| + |B^L|).
\]
So, if the negative term exceeds $2s^2$ then $G$ has fewer edges than the
construction $K_{s,a-s}\cup K_{s,b-s}$, a contradiction. Consequently we have $|A^M| + |B^M| + |B^L| < 4(s+1)^2$.
But then
$|E| \leq s|A^L| + (s-1/2)(|A^M| + |B^M| + |B^L|) < sa + 4(s+1)^3$
which is also a contradiction if $b$ is larger than $4(s+1)^3$. This completes the proof of (1).    

\medskip

As the unique extremal construction for (1) is not connected, we obtain $\exbc(a,b,D_{s,t})\le s(a+b-2s)-1$ if both $a$ and $b$ are at least $4(s+1)^3$ and $s\le t\le 2s-1$. On the other hand, one can alter the graph $K_{s,a-s}\cup K_{s,b-s}$ by removing two edges, one from each component, and adding an edge between the two vertices whose degree has been decreased to $s-1$ by the edge removal. As this new graph still has the property that no two vertices of degree at least $s+1$ are adjacent, it is $D_{s,t}$-free. The number of its edges is $s(a+b-2s)-2+1=s(a+b-2s)-1$.

\medskip

The proofs of (3), (4), and (5) use induction on $a+b$ and are analogous, so we show these statements together. First observe that \begin{itemize}
    \item 
    if $a\le s$, then $D_{s,t}\not\subseteq K_{a,b}$, and so $\exb(a,b,D_{s,t})=ab\le sb\le tb/2\le t(a+b)/2$;
    \item 
    if $a\le \ell$, then $S_{2\ell,(\ell+1) *1}\not\subseteq K_{a,b}$ and so $\exb(a,b,S_{2\ell,(\ell+1)*1})=ab\le \ell b\le \ell (a+b)$;
    \item 
    if $a<s+2$, then $P_{r,s,t}\not\subseteq K_{a,b}$ and so $\exb(a,b,P_{r,s,t})=ab\le m \cdot b\le m\cdot (a+b)$.
\end{itemize}  
So assume $a$ is large. Let $G\subseteq K_{a,b}$ be an $F$-free graph with $F$ being either $D_{s,t}$ or $S_{2(\ell+1),(\ell+1)*1}$ or $P_{r,s,t}$. If there exists a vertex $v$ of degree at most $s\le t/2$, $\ell$, $m$, respectively,  then by induction  we obtain that $|E(G)|$ is at most 
\begin{itemize}
    \item 
    $ t/2+E(G\setminus \{v\})\le t/2+\exb(a',b',D_{s,t})\le  t(a+b)/2$,
    \item 
    $\ell+E(G\setminus \{v\})\le \ell+\exb(a',b',S_{2\ell,(\ell+1)*1})\le  \ell (a+b)$,
    \item 
    $ m+E(G\setminus \{v\})\le m+\exb(a',b',P_{r,s,t})\le  m(a+b)$,
\end{itemize}
 as claimed, where either $a'=a$, $b'=b-1$ or $a'=a-1$, $b'=b$. 

On the other hand: 
\begin{itemize}
    \item 
    If all vertices have degree at least $s+1$, then $G$ cannot contain vertices of degree at least $t+1$, as then any such vertex and any of its neighbors would be centers of a copy of $D_{s,t}$. 
    \item 
    If all vertices have degree at least $\ell+1$, then we claim that $G$ cannot contain vertices of degree at least $2\ell+1$. Indeed, as $G$ is bipartite and all vertices have degree at least $\ell+1$, every vertex is an end vertex of a path $P_{2\ell+1}$ (we can select the vertices of $P_{2\ell+1}$ greedily). So if $v$ has degree $2\ell+1$, then we can choose $v$ as the end vertex of a path $P_{2\ell+1}$. Only $\ell$ vertices of the path can be neighbors of $v$ as $G$ is bipartite, so $v$ has at least $\ell+1$ other neighbors, giving a copy of $S_{2\ell,(\ell+1)*1}$.
    \item 
    If all vertices have degree at least $m+1$, then we claim that $G$ cannot contain vertices of degree at least $2m+1$. Indeed, suppose $m=t$ and $v$ has degree at least $2m+1$. Then let $N(v)=\{u_1,u_2,\dots,u_{2m+1},\dots\}$, $N(u_1)=\{v,v_1,\dots,v_{m},\dots\},$ and $N(v_1)=\{u_1,w_1,\dots,w_{m},\dots\}$. Then after picking a set $L_1$ of $r\le m$ neighbors of $v_1$ other than $u_1$, we can still pick a set $L$ of $m$ neighbors of $v$ with $L\cap (L_1\cup \{u_1\})=\emptyset$. As $G$ is bipartite, all neighbors of $u_1$ are outside $L_1\cup L$, and as $s<m$, we can pick $s$ of its neighbors that are distinct from $v,v_1$. In this way we obtain a copy of $P_{r,s,t}$ in $G$, a contradiction.
\end{itemize}

Therefore the degree sum is at most $t(a+b)$, $2\ell(a+b)$, and $2m(a+b)$, respectively, and so $|E(G)|$ is at most $t(a+b)/2$, $\ell(a+b)$, and $m(a+b)$. The sharpness of the bound follows from the lower bound of Proposition \ref{lower} (2).
\end{proof}

\begin{proof}[{\bf Proof of Theorem \ref{t:trees-k-l}}]
Part (1) is obtained as a corollary of Theorem \ref{classes} (1) and Proposition \ref{lower} (1), while Theorem \ref{classes} (3), (4) and Proposition \ref{lower} (1) yield part (2).
\end{proof}

\begin{proof}[{\bf Proof of Theorem \ref{gammab}}]
    To obtain the lower bound, it is enough to show that for any tree $T$ we have $\exb(n,T)\ge \frac{2}{3}(|T|-2)n-C$ where the constant $C$ depends only on $|T|$. Let the part sizes of $T$ be $k$ and $\ell$ with $k\le \ell$. Suppose first $\ell<2k$. Then by Proposition \ref{lower} (1) and $|T|-2=k+\ell-2\le 3k-3$, we have $\exb(n,T)\ge 2(k-1)(n-k+1)\ge \frac{2}{3}(|T|-2)n-|T|^2$. 
    
    If $2k\le \ell$, then by Proposition \ref{lower} (2) and $|T|-2\le \frac{3}{2}\ell-2$, we obtain $\exb(n,T)\ge \lfloor \frac{n}{\ell-1}\rfloor(\ell-1)^2\ge (\frac{2}{3}-\frac{1}{\ell})(|T|-2)n-C$, where $C$ depends only on $\ell$ and thus on $|T|$. This finishes the proof of the lower bound.

    The upper bound follows from any of Theorem \ref{classes} (1), (3), (4) or (5), if we pick $T=D_{k,2k-1}$, $T=D_{k,2k}$, $T=S_{2\ell,(\ell+1)*1}$ or $T=P_{k,k-1,k}$. Note that in all cases the larger part in the bipartition of $T$ is roughly twice the size of the smaller part.
\end{proof}

\begin{proof}[{\bf Proof of Theorem \ref{gammabc}}] Let us fix  constants $c$ and $p$ with $2/3 \le c<1$ and $p\le c$. Consider the following three constructions:
\begin{enumerate}
    \item 
    Let $K^1,K^2,\dots$ be copies of $K_{(2c-1)k,(2c-1)k}$, and let $P^1,P^2,\dots$ be copies of $P_{2(1-c)k}$. Let $G$ be the connected subgraph of $K_{n,n}$ with parts $A$ and $B$ such that the subgraph induced on the first $2(c-1)k$ vertices on both parts is $K^1$, the subgraph induced on the next $(1-c)k$ vertices on both sides is $P^1$ with the first vertex being in $A$ and last vertex in $B$, and then $K^i$, $P^i$ alternate. To achieve connectivity, we join a vertex of $K^i$ in $B$ to the start vertex of $P^i$ in $A$, and the endvertex of $P^i$ in $B$ to a vertex of $K^{i+1}$ in $A$. We claim that $G$ does not contain any tree $T$ on $k$ vertices with the larger part of $T$ containing more than $ck$ vertices. Indeed, a copy of $T$ in $G$ cannot contain vertices from $K^i$ and $K^{i+1}$ as then it should contain all of $P^i$ that contains $(1-c)k$ vertices in both parts. On the other hand if a copy of $T$ contains vertices only from $P^{i-1}\cup K^i\cup P^i$, then to have $ck$ vertices on one side, it needs to contain at least $ck-(2c-1)k=(1-c)k$ vertices from $P^{i-1}\cup P^i$, but these vertices come in the same amount on both parts, so the copy would contain $(1-c)k$ vertices on the other side contradicting that the smaller size is strictly smaller than $(1-c)k$. The number of edges in $G$ is at least $(2c-1)^2k^2\cdot \lfloor \frac{n}{(2c-1)k+(1-c)k}\rfloor\ge \frac{(2c-1)^2k}{c}n$, an increasing function of $c$ in the range $c>1/2$.
    \item 
    Suppose $T$ contains a path longer than $pk$. Then let $G$ be the subgraph of $K_{n,n}$ with parts $A$ and $B$ containing universal vertices (vertices joined to all vertices of the other part) $u\in A$ and $v_2,v_3,\dots,v_{\lfloor \frac{pk}{2}\rfloor}\in B$ and no other edges. As the longest path in $G$ contains $2\lfloor \frac{pk}{2}\rfloor\le pk$ vertices, $G$ is $T$-free and the number of edges in $G$ is $(\lfloor \frac{p}{2}k\rfloor -o(1))n$. 
    \item 
    Suppose finally that $T$ does not contain a path of  more than $pk$ vertices and the larger part of $T$ contains at most $ck$ vertices. Let $\ell$ be the smallest even number larger than $pk$ and consider the subgraph $G\subseteq K_{n,n}$ with parts $A$ and $B$ that is defined as follows: partition $A$ into $A_1,A_2,A_3$, $B$ into $B_1,B_2,B_3$ with $|A_1|=|B_1|=\frac{\ell}{2}$ and $|A_2|=|B_2|=\alpha=(1-c-\frac{p}{2})k-2$, and let $G_1=G[A_1\cup B_1]=P$ a path on $\ell$ vertices with endvertices $u\in A$ and $v\in B$. Further, let $G_2=G[A_2,B_3\cup \{v\}]$ and $G_3=G[A_3\cup \{u\},B_2]$ be complete bipartite graphs and the graphs $G_1,G_2,G_3$ cover all edges of $G$.
    We claim that $G$ is $T$-free. First observe that as $T$ does not contain a path on $\ell$ vertices, therefore a copy of $T$ in $G$ cannot contain vertices from both components of $G\setminus V(P)$. On the other hand, if $C$ denotes one of the components of $G\setminus V(P)$, the smaller part of $G[C\cup V(P)]$ has size at most $\alpha+ \frac{\ell}{2}<(1-c)k$, which is less than the size of the smaller part of $T$, so $G$ is $T$-free as claimed. The number of edges in $G$ is $2(1-c-\frac{p}{2}-o(1))kn$.
\end{enumerate}
    The rest of the proof is a maximization problem. Given $c$, the function $\min\{\frac{p}{2},2(1-c-\frac{p}{2})\}$ is maximized when $p=\frac{4-4c}{3}$, with value $\frac{2-2c}{3}$. As $\frac{(2c-1)^2}{c}$ is an increasing function for $c>1/2$, the value $\min_c\max\{\frac{2-2c}{3},\frac{(2c-1)^2}{c}\}$ is attained when the two expressions are equal, yielding the quadratic  equation $14x^2-14x+3=0$. So if $x_0$ is the larger root of this equation, then setting $c=x_0$ and $p=\frac{4-4x_0}{3}$ one obtains that 
    \begin{itemize}
        \item 
        $\frac{\exbc(n,T)}{(|T|-2)n}$ is at least $\frac{(2x_0-1)^2}{x_0}-o(1)$ via the first construction if $T$'s larger part contains more than $x_0k$ vertices,
        \item 
        $\frac{\exbc(n,T)}{(|T|-2)n}$ is at least $\frac{2-2x_0}{3}-o(1)$ via the second construction if $T$ contains a path with more than $\frac{4x_0-4}{3}$ vertices,
        \item 
        $\frac{\exbc(n,T)}{(|T|-2)n}$ is at least $2(1-x_0-\frac{2-2x_0}{3})-o(1)$ via the third construction, if $T$ does not  fall into the above categories, 
    \end{itemize}
    where the $o(1)$ term tends to 0 as $|T|$ tends to infinity. As these three values are equal, this finishes the proof of the theorem.
\end{proof}

\bsk

\section{Proofs for small trees}
\label{s:small}

\begin{proof}[{\bf Proof of Theorem \ref{smalltrees}}]
The upper bound in (1) is the particular case $k=2$ of Theorem \ref{spiders} (2), and tightness for $3\mid n$ is shown by $\frac{n}{3}K_{3,3}$.
No other extremal graphs exist if $n$ is a multiple of~3, as it follows from the proof of Theorem \ref{spiders}.

\medskip

The assertion of (2) is the particular case $a=b=n$ and $s=2$ of Theorem \ref{classes} (1), provided that $n$ is not too small (more explicitly $n\geq 108$).
For smaller $n$ a careful analysis is needed.
Before that, we note that $2K_{2,n-2}$ is a universal construction yielding the lower bound $4n-8$ for all $n$.
However, the graph $\theta_{(n-1)*3}$ with two nonadjacent vertices $u,v$ of degree $n-1$ and a perfect matching connecting $N(u)$ with $N(v)$ is $D_{2,2}$-free and has $3n-3$ edges.
Hence, it is an alternative extremal graph if $n=5$ and has more than $4n-8$ edges if $n=3,4$.
Further, $C_8$ has $8=4n-8$ edges for $n=4$ (fewer than $3n-3$) and $C_6$ has $6>4n-8$ edges for $n=3$ (same as $3n-3$, in fact $C_6\cong \theta_{2*3}$).

Let now $G \subset K_{n,n}$ be a $D_{2,2}$-free bipartite graph with vertex classes $A,B$.
We write $A=A^+\cup A^-$ and $B=B^+\cup B^-$, where $A^+\cup B^+$ is the set of vertices of degree at least 3.
Any edge joining $A^+$ with $B^+$ would create a copy of $D_{2,2}$, hence $A^+\cup B^+$ is an independent set.
It follows that $e(G)\leq 2|A^-|+2|B^-|$, which yields at most $4n-10$ edges whenever $|A^+\cup B^+|\geq 5$. Also, $e(G) = \sum_{v\in A} d(v) \leq 2n$ if $A^+=\emptyset$ (or analogously if $B^+=\emptyset$) because $A$ (and also $B$) covers all edges.
Then $2n\leq 4n-8$ for $n\geq 4$, and the computation leads to the cycle $C_6$ and the disconnected graph $K_{2,3} \cup K_1$ (both are extremal) with $6$ if $n=3$.

Consider next $|A^+|=1$, $A^+=\{v'\}$.
Then $e(G) = d(v') + \sum_{v\in A\setminus \{v'\}} d(v) \leq |B^-| + 2 |A^-| \leq 3n-3$, with equality only if $|B^-|=n-1$ and all $v\in A\setminus \{v'\}$ have degree 2.
As $|B^+|=1$, also all $v''\in B^-$ have degree 2.
Consequently $G\cong \theta_{(n-1)*3}$ in this case.
Its number of edges exceeds $4n-8$ only if $n\leq 4$.

The only remaining case is $|A^+| = |B^+| = 2$.
Since $A^+\cup B^+$ is independent, we have $e(G)\leq 2|A^-| + 2|B^-| \leq 4n-8$, with equality only if also $A^-\cup B^-$ is an independent set.
Thus, $G\cong 2K_{2,n-2}$, completing the proof of (2).

\medskip

Concerning $S_{2,2,1}$ suppose that $G\subset K_{n,n}$ is a graph with $4n-7$ edges, $n\geq 4$.
If $G\cong\theta_{3*3}$, then $S_{2,2,1}\subset G$.
Otherwise $D_{2,2}\subset G$ by the above proof of (2).
Let the central edge of this double star be $uv$, and its leaves $u_1,u_2,v_1,v_2$.
If any of these leaves has degree bigger than 1, we directly find a copy of $S_{2,2,1}$, because $G$ is bipartite.
Otherwise let $G'=G-u_1-u_2-v_1-v_2$, and write $n'=n-2$.
Now we have $G'\subset K_{n',n'}$ and $e(G')=4n'-3$.
We cannot have $n'=2$ (i.e., $n=4$) because then $4n'-3=5>e(K_{2,2})$.
If $n'=3$ (i.e., $n=5$), then $4n'-3=9=e(K_{3,3})$, and obviously $S_{2,2,1}\subset K_{3,3}$.
As a consequence, the inequality $\exb(n,S_{2,2,1})\leq 4n-8$ is valid for $n=4,5$.
From this basis we can prove (3) by the method of infinite descent.
If (3) is false, let $n$ be the smallest integer for which there exists an $S_{2,2,1}$-free subgraph of $K_{n,n}$ with $4n-7$ edges.
On applying the above reduction, we obtain the contradiction that $\exb(n',S_{2,2,1})\geq 4n'-3$ holds for $n-2=n'\geq 4$.

Finally, in case of $n=3$, deleting any two edges from $K_{3,3}$, the remaining graph obviously contains $S_{2,2,1}$, hence $\exb(3,S_{2,2,1})\leq 6$; and this is tight due to the fact $S_{2,2,1}\not\subset C_6$.
\end{proof}

\begin{proof}[{\bf Proof of Theorem \ref{smalltreesconnected}}]
The lower bound $2n-1$ in (1) is shown by the double star $D_{n-1,n-1}$.
To see the corresponding upper bounds, observe that a graph $G$ on $2n$ vertices with $2n$ edges contains a cycle $C$. If $C$ has length at least $6$, then $G$ contains a $P_6$. If $C$ is shorter, then $C$ has length 4 as $G$ is bipartite, and thus by connectivity, if there exist vertices outside $C$ on both sides, which is the case if $n\ge 3$, $G$ contains a $P_6$.

\medskip

Concerning (2), we have seen in the proof of Theorem
\ref{smalltrees} that the graph $\theta_{(n-1)*3}$ is extremal for $\exb(n,D_{2,2})$ if $3\leq n\leq 5$.
Since it is a connected graph, it also correctly gives the value of $\exbc(n,D_{2,2})$.
On the other hand, if $n\geq 6$, then the unique extremal graph for $\exb(n,D_{2,2})$ is disconnected, hence $\exbc(n,D_{2,2})\leq 4n-9$ follows.
This upper bound is tight for all $n\geq 6$, as shown by the $D_{2,2}$-free construction on $u_1,u_2,\dots,u_n$, $v_1,v_2,\dots,v_n$ with edges $u_1v_1,u_1v_2,u_2v_1$ and $u_2v_j,u_3v_j,u_jv_2,u_jv_3$ for all $4\le j\le n$.

\medskip

The lower bound of (3) is given by a connected 2-factor (i.e., Hamiltonian cycle) of $K_{n,n}$.
For the upper bound of (3) suppose $G\subseteq K_{n,n}$ is connected  with $e(G) \ge 2n +1$. Then $G$ contains an even cycle.  If the cycle is a spanning $C_{2n}$, then it must have a chord as $e(G)>2n$ and this chord produces a copy of $S_{3,1,1}$ if $n \ge 4$. 
(Note that $C_6$ and an antipodal chord does not contain $S_{3,1,1}$.)
If the cycle $C$ is of length $2k$ with $6 \le 2k \le 2n -2$, then by connectivity of $G$ at least one vertex $v\in C$  must have a neighbor $u\notin C$ giving a copy of $S_{3,1,1}$ with center at $v$.

Suppose that $G$ contains only cycles of length 4. As $e(G) \ge 2n +1$, there must exist two such cycles. Suppose first the two $C_4$s are vertex-disjoint.  Then  by connectivity there is a path  connecting them, creating a copy of $S_{3,1,1}$. If the two $C_4$s  have one vertex in common, then that vertex is the center of a copy of $S_{3,1,1}$. If the two $C_4$s have two vertices in common, but not an edge, then they form a $K_{2,4}$ and there is $S_{3,1,1}$ with center at a vertex of degree~4.  
And if the two $C_4$s have two vertices in common that are adjacent, then their union contains $C_6$, a case that we have already dealt with. 

Finally, if the two $C_4$s have three common vertices, they form a copy of $K_{2,3}$  with parts say $\{ x,y \}$  and $\{ a,b,c\}$. If one of the remaining vertices is adjacent to $x$ or $y$, then there is a copy of $S_{3,1,1}$ centered at $x$ or $y$.
If vertex $a$ has two  neighbors apart from $x,y$, then there is $S_{3,1,1}$ centered at $a$, and the same holds for $b$ and $c$.  If $a$ and $b$ have distinct neighbors different from $x$ and $y$, then there is a copy of $S_{3,1,1}$ and the same holds for $a ,c$ and $b ,c$. Hence $a, b, c$  can be adjacent to at most one vertex $w$ apart from $x$ and $y$.

As $n \ge 4$, there is a further vertex $z$ in the side of $x$ and $y$. If $x$ or $y$ is an end vertex of a $P_4$ using at most one of $a,b,c$, then $G$ contains a $S_{3,1,1}$. If $w$ does not exist, then without such a path, the side of $x$ and $y$ is smaller than the other side. If $w$ does exist, then $zdx$ or $zdy$ or $zdw$ is a path and all these would create an $S_{3,1,1}$. This finishes the upper bound of (3).

\medskip

The lower bound of (4) is again a connected 2-factor of $K_{n,n}$. The proof of the upper bound of (4) is very similar to that of (3). Suppose $G\subseteq K_{n,n}$ has $2n+1$ edges. Any even cycle of length at least 6 with an extra chord or a pendant edge contains $S_{2,2,1}$, so we can assume that $G$ contains only $C_4$s, and the number of edges ensures that there should be at least two such $C_4$s.

If the two $C_4$s are vertex disjoint, then  by connectivity there is a path connecting them  and creating an $S_{2,2,1}$. If the two $C_4$s have one vertex in common, then it is the center of an $S_{2,2,1}$. If the two $C_4$s have two vertices in common, but not an edge, their union forms $K_{2,4}$ with $A = \{ x,y \}$ vertices of degree 4 and $B =\{  a,b,c,d\}$ vertices of degree 2.

Consider a vertex $w$ not in the $K_{2,4}$ but adjacent to some vertex of $B$. Then there is an $S_{2,2,1}$.   So $a,b,c,d$ are only connected to $x$ and $y$. Consider $w$ in the $A$-side. By connectivity, it must have a path of length at least 2 to some vertex of $A$, say to $x$. Then $x$ is the center of $S_{2,2,1}$ settling the case of edge-disjoint $C_4$s.

If the two $C_4$s have exactly two adjacent vertices in common, then they form $C_6$ and an antipodal chord, containing $S_{2,2,1}$.   

If the two $C_4$s have three common vertices, they form a $K_{2,3}$, say with $A$-part  $\{ x,y \}$  and $B$-part $\{ a,b,c\}$.  If one of the remaining vertices is adjacent to a member of $B$, there is $S_{2,2,1}$. So consider a vertex $w$ in the $A$-side. By connectivity it can reach the $K_{2,3}$ via $A$ hence there is a path of length at least 2 from $w$ to some vertex of $A$, say to $x$. Then $x$ is the center of an $S_{2,2,1}$.
\end{proof}

\section{Concluding remarks and open problems}

We have started a systematic study of $\exb(a,b,F)$ and $\exbc(a, b, F)$ giving 
general lower-bound constructions, and getting exact values of these functions 
for small trees up to six vertices as well as in many cases of double stars 
and spiders and the family of trees $\cT_{k,\ell}$  proposed in \cite{YZ}.

We also considered (as done in \cite{CPT} and \cite{JLS}) worst-case scenarios solving 
$\gamma_b = 2/3$ and getting lower bounds for $\gamma_{b,c}$ and $\gamma_{b\leftrightarrow c}$.

Perhaps the most important open problems emerging from this paper, in this initial 
stage, are:

\begin{problem}
Determine exactly or improve upon the lower bounds obtained for
$\gamma_{b,c}$ and $\gamma_{b\leftrightarrow c}$.
\end{problem}

\begin{problem}
Determine $\exb(a,b,F)$ and $\exbc(a, b, F)$ for further families of 
trees, in particular for spiders of the form $S_{k,d*1}$ extending the already solved 
cases in Theorems \ref{spiders} and \ref{classes}.
\end{problem}

\begin{problem}
Determine $\exb(a,b,F)$ and $\exbc(a, b, F)$ for all the trees on 
7 vertices,  with the hope that some of the proofs extend to infinite families 
of trees (as was the case that led to Theorems \ref{spiders}, \ref{classes}, and \ref{t:trees-k-l}).
\end{problem}

\begin{problem}
    Determine at least asymptotically $\exbc(n,P_k)$ for any fixed $k$. This would be an extension of Theorem \ref{smalltreesconnected} (1) that establishes $\exbc(n,P_5)=\exbc(n,P_6)=2n-1$.
\end{problem}

\end{document}